\documentclass[final,3p,times]{elsarticle}
\usepackage{}

\usepackage{graphicx,amsmath,amssymb,txfonts}
\usepackage{algorithm}
\usepackage{algorithmic}
\usepackage{subfigure,color}
\usepackage[colorlinks]{hyperref}
\usepackage{hypernat}
\usepackage{multirow}
\numberwithin{equation}{section}
\graphicspath{{Figs/}}
\usepackage{booktabs}
\usepackage{listings}
\usepackage{arydshln}
\usepackage{epstopdf}

\biboptions{sort&compress}

\makeatletter

\newcommand{\Rmnum}[1]{\expandafter\@slowromancap\romannumeral #1@}
\makeatother

\journal{Global Science Press}

\newtheorem{definition}{\textbf{Definition}}

\newtheorem{remark}{\textbf{Remark}}

\newcommand{\ep}{\hfill\rule{0.15cm}{0.35cm}\vskip 0.3cm}

\begin{document}

\begin{frontmatter}



\title{Algebraic multigrid block preconditioning for 
multi-group radiation diffusion equations}


\author[myadda]{Xiaoqiang Yue}
\ead{yuexq@xtu.edu.cn}

\author[myadda]{Shulei Zhang}

\author[myaddb]{Xiaowen Xu}
\ead{xwxu@iapcm.ac.cn}

\author[myadda]{Shi Shu\corref{mycoraur}}
\ead{shushi@xtu.edu.cn}



\author[myaddd]{Weidong Shi}

\address[myadda]{School of Mathematics and Computational Science,
Hunan Key Laboratory for Computation and Simulation in Science and Engineering,
Key Laboratory of Intelligent Computing \& Information Processing of Ministry of Education,
Xiangtan University, Xiangtan 411105, China}
\address[myaddb]{Laboratory of Computational Physics,
Institute of Applied Physics and Computational Mathematics,
Beijing 100088, China}
\address[myaddd]{School of Applied Mathematics,
Shanxi University of Finance and Economics, Taiyuan 030006, China}

\cortext[mycoraur]{Corresponding author.}

\begin{abstract}
The paper focuses on developing and studying efficient block preconditioners based on classical algebraic multigrid
for the large-scale sparse linear systems
arising from the fully coupled and implicitly cell-centered finite volume discretization of
multi-group radiation diffusion equations, 
whose coefficient matrices can be rearranged into the $(G+2)\times(G+2)$ block form, where $G$ is the number of energy groups.
The preconditioning techniques are based on the monolithic classical algebraic multigrid method,
physical-variable based coarsening two-level algorithm and two types of block Schur complement preconditioners.
The classical algebraic multigrid is applied to solve the subsystems that arise in the last three block preconditioners.
The coupling strength and diagonal dominance are further explored to improve performance.
We use representative one-group and twenty-group linear systems from capsule implosion simulations to test
the robustness, efficiency, strong and weak parallel scaling properties of the proposed methods.
Numerical results demonstrate that block preconditioners lead to mesh- and problem-independent convergence,
and scale well both algorithmically and in parallel.
\end{abstract}


\end{frontmatter}

\section{Introduction}\label{sec1}

The multi-group radiation diffusion (MGD) equations have a broad range applications,
including inertial confinement fusion (ICF) and astrophysics \cite{m-04}.
As a result of the complicated nonlinear couplings among dozens of physical quantities at multiple temporal and spatial scales,
MGD equations are often discretized by the finite volume method allowing for local conservations \cite{h-03,n-01,h-04,g-01,p-03,s-02,x-05},
resulting in a series of nonsymmetric but positive definite large-scale sparse linear systems.
The overall complexity increases not only with mesh sizes but also with the level of couplings between physical quantities.
It must be emphasized that these numerical solutions play a time-consuming role (about 80\% in general) in ICF numerical simulations,
due to the fact that the coefficient matrices are invariably ill-conditioned.

To effectively address the aforesaid bottlenecks, numerous preconditioned Krylov subspace methods
have been proposed in an efficient and scalable manner over the past decades,
see \cite{b-01,m-01,x-03,y-01,y-02,x-01,x-02,h-02,y-03} and references cited therein.
These preconditioners are conceived as approximate inverses and they
fall into the category of incomplete LU factorizations,
domain decomposition preconditioners, monolithic algebraic multigrid (AMG) methods
and their symmetric / nonsymmetric combinations.
Since each of those coefficient matrices has an underlying block structure,
one can also determine a block preconditioner to separate the global problem into easier-to-solve subproblems
and form an object-oriented framework to allow for code-reuse and
incorporate single-physics experience into multi-physics simulations.
The objective is to approximately invert numerous individual scalar systems instead of the fully coupled systems.
Preconditioners of this type had been proposed and analyzed in the literature,
including the block diagonal preconditioner \cite{d-01,l-01,z-03},
block lower / upper triangular preconditioner \cite{b-04,b-02,c-01},
product (splitting) preconditioner \cite{r-02,w-01,z-02} and constraint preconditioner \cite{b-06,k-01,d-02}.
Block preconditioners with multigrid components had proven very successful in a variety of applications,
e.g., liquid crystal directors modeling \cite{b-07},
multiphase flow in porous media \cite{b-03}, Stokes problem \cite{c-02},
incompressible Navier-Stokes problem \cite{e-01},
second-order Agmon-Douglis-Nirenberg elliptic systems \cite{l-02},
magnetohydrodynamics model \cite{m-02},
Dirichlet biharmonic problem \cite{p-02},
electrical activity in the heart \cite{s-01},
Brinkman problem \cite{v-01}, all-speed melt pool flow physics \cite{w-03} and fully coupled flow and geomechanics \cite{w-02}.
Our focus in this work is on the block preconditioning
based on the classical AMG method owing to its general applicability, high efficiency and easy-to-use user interface.

In the recent work \cite{a-01}, four types of operator-based preconditioners
have been developed in the Jacobian-free Newton-Krylov method
for solving two-dimensional three-temperature energy equations.
These preconditioners are application specific,
relying on physical properties of the energy equations as well as
different linearizations on different terms in the nonlinear residual.
They are demonstrated numerically to be very effective. However,
the corresponding preconditioning matrix has to be assembled explicitly in each time step.
In this study, we restrict our consideration to a variety of purely algebraic block preconditioners,
which are based only on information available from the original block system,
without constructing the preconditioning matrix.

The arrangement of this work proceeds as follows. Section \ref{sec2}
provides the mathematical formulation for MGD equations and the fully coupled and implicitly cell-centered finite volume discretization.
Section \ref{sec3} describes several AMG block-based preconditioning strategies,
including an adaptive preconditioning strategy and four approximations on matrix inverse
to improve performance. In Section \ref{sec4}, we compare the performance of these preconditioners for
representative MGD linear systems from capsule implosion simulations, on both sequential and parallel computers.
Finally, we discuss conclusions in Section \ref{sec5}.

\section{Problem formulation and discretization}\label{sec2}

This article is concerned with the time-dependent MGD equations in a certain symmetric geometry \cite{p-01}
\begin{align}\label{eq-2-01}
  \left \{
    \begin{aligned}
      & \frac{\partial E_g}{\partial t} = \nabla \cdot (D_g(E_g)\nabla E_g) +
    c(\sigma_{Bg}B(T_E)_g - \sigma_{Pg}E_g) + S_g,\qquad g=1,\cdots,G \\
      & \rho c_E\frac{\partial T_E}{\partial t} = \nabla \cdot (D_E(T_E)\nabla T_E) -
    c\sum_{g=1}^G(\sigma_{Bg}B(T_E)_g - \sigma_{Pg}E_g) + \omega_{IE}(T_I - T_E) \\
      & \rho c_I\frac{\partial T_I}{\partial t} = \nabla \cdot (D_I(T_I)\nabla T_I) - \omega_{IE}(T_I - T_E)
    \end{aligned}
  \right.,
\end{align}
where
\begin{itemize}

  \item $G$, $c$, $\rho$, $\omega_{IE}$ respectively denote the number of energy groups, velocity of light,
  density of the medium and energy transfer coefficient between electron and ion;

  \item $E_g$, $B(T_E)_g$, $D_g(E_g)$ respectively denote the $g$-th spectral radiation
  and electron scattering energy densities and highly nonlinear radiation diffusion coefficient;

  \item $S_g$, $\sigma_{Bg}$, $\sigma_{Pg}$ respectively denote the $g$-th radiation source,
  scattering and absorption coefficients of the Planck-averaged electron energy;

  \item $c_\alpha$, $T_\alpha$, $D_\alpha(T_\alpha)$ respectively denote the specific heat capacity,
  temperature and nonlinear thermal-conductivity coefficient of electron ($\alpha$ = $E$) or ion ($\alpha$ = $I$).

\end{itemize}

Various discretization schemes (see \cite{q-01}) can be applied to reduce
the continuous differential equations to finite dimensional sparse linear systems.
Utilizing the (adaptive) backward Euler for the temporal discretization,
followed by the frozen-in coefficients for the linearization,
and then an appropriate cell-centered finite volume for the spatial discretization,
we obtain a series of sparse $(G+2)\times(G+2)$ block structured linear systems from \eqref{eq-2-01}
by grouping together the unknowns corresponding to the same physical quantity:
\begin{eqnarray}\label{eq-2-02}
  {\bf A} {\bf T} \equiv
  \left[
  \begin{array}{ccccc}
    A_1    &        &        & D_{1E} &  \\
           & \ddots &        & \vdots &  \\
           &        & A_G    & D_{GE} &  \\ 
    D_{E1} & \cdots & D_{EG} & A_E & D_{EI} \\ 
           &        &        & D_{IE} & A_I
  \end{array}\right]\left[
  \begin{array}{c}
    T_1    \\
    \vdots \\
    T_G    \\ 
    T_E    \\ 
    T_I
  \end{array}\right] = \left[
  \begin{array}{c}
    f_1    \\
    \vdots \\
    f_G    \\ 
    f_E    \\ 
    f_I
  \end{array}\right]
  \equiv {\bf f},
\end{eqnarray}
in which $D_{EI} = D_{IE}$ and $D_{Eg} \ne D_{gE}$ for group index $g=1,\cdots,G$,
causing the coefficient matrix is generally positive definite but necessarily nonsymmetric.
For these linear systems, we use the restarted generalized minimal residual (GMRES($m$)) method,
where $m$ is the number of Krylov directions to orthogonalize against.
Furthermore, preconditioning the GMRES($m$) solver is essential for rapid convergence.

It is worth noting that the diagonal blocks $A_\alpha$ ($\alpha=1,\cdots,G,E,I$)
have the same nonzero structure of a discrete purely elliptic problem,
while the coupling terms $D_{\alpha\beta}$ ($\alpha\ne\beta$) are nonsingular diagonal matrices.
Typically, algebraic characteristics of submatrices $A_\alpha$ ($\alpha=1,\cdots,G,E,I$)
are much better than those of the matrix ${\bf A}$.
This is the fundamental motivation for devising block preconditioners to solve \eqref{eq-2-02}.

\section{Preconditioning strategies}\label{sec3}

This section is devoted to four preconditioning strategies:
the monolithic classical AMG, physical-variable based coarsening two-level (PCTL)
and two types of block Schur complement preconditioners.
Furthermore, two further improvements are introduced in Sections \ref{sec3.5} and \ref{sec3.6}.

\subsection{Monolithic classical AMG preconditioner}\label{sec3.1}

Nowadays classical AMG developed in \cite{r-01} is quite mature and one of the most popular preconditioners in real applications,
since its virtue of scalability and applicability on complicated domains and unstructured grids. It has Setup and Solve phases.
The former phase builds all the ingredients required by a hierarchy of grids, the finest to the coarsest,
under the assumption that no information on the underlying geometry, grids and continuous operators is available,
while the latter phase performs V-cycle, W-cycle or F-cycle.
Regarding more detailed information and challenges, we refer to the review articles \cite{s-04,x-07,x-06} and related references therein.
It should be noted, however, that different coarsening schemes, interpolation procedures and relaxation methods
result in different operator complexities, Setup time, cycle time and convergence rates.
We apply BoomerAMG \cite{h-01}, a parallel implementation of classical AMG in the HYPRE package,
on the fully-coupled monolithic linear system as a black-box preconditioner,
empirically with a strength threshold 0.25, 
a single V-cycle with Falgout (a hybrid RS / CLJP) \cite{h-01} 
coarsening for 
test runs, an aggressive coarsening process on the finest level,
one presmoothing and postsmoothing sweep performed by the hybrid Gauss-Seidel in symmetric ordering
(i.e., down cycle with relaxations swept first all coarse points and then all fine points,
and up cycle using the opposite traversal), classical modified interpolation,
at most 100 degrees of freedom (DoFs) at the coarsest level which is solved via Gaussian elimination,
and other parameters chosen from the default configuration.

\subsection{PCTL preconditioner}\label{sec3.3}

The PCTL algorithm was first proposed by Xu and his co-authors for 3-temperature ($G=1$) linear systems \cite{x-03}.
The four-level multigrid reduction preconditioner proposed recently by Bui, Wang and Osei-Kuffuor \cite{b-08}
can be viewed as a variation of PCTL. The generalization of PCTL used to precondition \eqref{eq-2-02} is straightforward.
Various components have to be chosen.
Taking the special structure of the coefficient matrix into account,
we set the electron temperature variables to be coarse points,
and the others as fine points. The relaxation routine used here is C / F block relaxation sweeps, i.e.,
variables associated with coarse points are relaxed followed by the other unknowns.
Denote by
\begin{eqnarray*}
 {\bf P} = \left[
  \begin{array}{ccccc}
    P_1^\top & \cdots & P_G^\top & I & P_I^\top
  \end{array}\right]^\top
\end{eqnarray*}
the interpolation matrix to match up with the aforementioned coarse-grid selection,
where submatrices $P_g = -A_g^{-1}D_{gE}$ ($g=1,\cdots,G$) and $P_I = -A_I^{-1}D_{IE}$.
However, they are often significantly denser, even though $D_{\alpha E}$ are all diagonal matrices. 
From the implementation point of view, they are restricted to be diagonal
and give an exact approximation only for constant functions.
That is, for vector ${\bf 1}=(1,\cdots,1)^\top$,
the resulting vectors containing all the diagonal entries of blocks $P_\alpha$ ($\alpha=1,\cdots,G,I$),
denoted by $p_\alpha$, can be respectively computed as follows:
\begin{eqnarray*}
 p_g = -A_g^{-1}D_{gE}{\bf 1},~g=1,\cdots,G;\quad p_I = -A_I^{-1}D_{IE}{\bf 1}.
\end{eqnarray*}
The restriction is defined as ${\bf P}^\top$ and the coarse-grid operator $A_c$ is taken as
\begin{eqnarray}\label{eq-3-03}
  A_c = {\bf P}^\top {\bf A}{\bf P} = \sum_{g=1}^G P_g^\top A_g P_g + A_E
  + P_I^\top A_I P_I + \sum_{g=1}^G D_{Eg} P_g + \sum_{g=1}^G P_g^\top D_{gE} + D_{EI} P_I + P_I^\top D_{IE}.
\end{eqnarray}
Noting that $P_\alpha$ ($\alpha=1,\cdots,G,I$) are all diagonal,
it can be deduced that $A_c$ and $A_E$ have the same nonzero structure.
The preconditioning phase of the PCTL algorithm takes the following form:
\begin{enumerate}

  \item (pre-smoothing step) Run block relaxation once to subsystems:
\begin{eqnarray*}
  w_E = A_E^{-1} b_E;\quad w_g = A_g^{-1} (b_g - D_{gE}w_E),~g=1,\cdots,G;\quad w_I = A_I^{-1} (b_I - D_{IE}w_E);
\end{eqnarray*}

  \item Solve the coarse-grid system: $A_cw_c = {\bf P}^\top ({\bf b} - {\bf A} {\bf w})$;

  \item (coarse-grid correction) Set ${\bf w}:={\bf w}+{\bf P}w_c$.

\end{enumerate}
where ${\bf b}=(b_1,\cdots,b_G,b_E,b_I)^\top$ is an incoming arbitrary Krylov vector and
${\bf w}=(w_1,\cdots,w_G,w_E,w_I)^\top$ is the outgoing Krylov solution.

\subsection{Two types of block Schur complement preconditioners}\label{sec3.4}

Noting that the coefficient matrix ${\bf A}$ can be rewritten in a factored form
\begin{eqnarray}\label{eq-3-01}
{\bf A} = \left[
  \begin{array}{ccccc}
    I    &        &        &  &  \\
           & \ddots &        &  &  \\
           &        & I    &  &  \\ 
     &  &  & I & D_{EI}A_I^{-1} \\ 
           &        &        &  & I
  \end{array}\right]\left[
  \begin{array}{ccccc}
    A_1    &        &        & D_{1E} &  \\
           & \ddots &        & \vdots &  \\
           &        & A_G    & D_{GE} &  \\ 
    D_{E1} & \cdots & D_{EG} & C_E &  \\ 
           &        &        &  & A_I
  \end{array}\right] \left[
  \begin{array}{ccccc}
    I    &        &        &  &  \\
           & \ddots &        &  &  \\
           &        & I    &  &  \\ 
     &  &  & I &  \\ 
           &        &        & A_I^{-1} D_{IE} & I
  \end{array}\right],
\end{eqnarray}
where the matrix $C_E = A_E - D_{EI}A_I^{-1}D_{IE}$ and $I$ is the identity matrix with the same size of $A_I$.
At this point, it is important to emphasize that the factorization \eqref{eq-3-01}
breaks $A_I$ away from the original system. Then a similar factorization
\begin{eqnarray*}
  \left[
  \begin{array}{cccc}
    A_1    &        &        & D_{1E} \\
           & \ddots &        & \vdots \\
           &        & A_G    & D_{GE} \\
    D_{E1} & \cdots & D_{EG} & C_E
  \end{array}\right]=\left[
  \begin{array}{cccc}
    I    &        &        & D_{1E}C_E^{-1} \\
           & \ddots &        & \vdots \\
           &        & I    & D_{GE}C_E^{-1} \\
     &  &  & I
  \end{array}\right]\left[
  \begin{array}{ccc;{0.5pt/2pt}c}
        &        &        &  \\
           & C_R &        &  \\
           &        &     &  \\ \hdashline[0.5pt/2pt]
     &  &  & C_E
  \end{array}\right]\left[
  \begin{array}{cccc}
    I    &        &        &  \\
           & \ddots &        &  \\
           &        & I    &  \\
    C_E^{-1}D_{E1} & \cdots & C_E^{-1}D_{EG} & I
  \end{array}\right]
\end{eqnarray*}
is processed to decouple the remaining in two separate systems,
where
\begin{eqnarray*}
  C_R = \left[ \begin{array}{ccc}
    A_1    &        &    \\
           & \ddots &    \\
           &        & A_G
  \end{array}\right]-\left[
  \begin{array}{c}
    D_{1E} \\
    \vdots \\
    D_{GE}
  \end{array}\right] C_E^{-1} \left[
  \begin{array}{ccc}
    D_{E1} & \cdots & D_{EG}
  \end{array}\right].
\end{eqnarray*}
The error associated with the global factorization would be zero if $C_E$ is inverted exactly
and $C_R$ is determined completely. However, it would be prohibitively expensive to explicitly form the full block matrix $C_R$,
we instead implement $\tilde{C}_R = {\bf diag}(\tilde{C}_1,\cdots,\tilde{C}_G)$
with $\tilde{C}_g=A_g-D_{gE}C_E^{-1}D_{Eg}$ for $g=1,\cdots,G$,
ignoring all off-diagonal elements of $C_R$.
In this situation, we obtain the completely decoupled shape,
leading to the following specific implementation procedure:
\begin{enumerate}

  \item (the intermediate ion segment) $w_I^* = A_I^{-1}b_I$;

  \item (the intermediate electron segment) $w_E^* = C_E^{-1}(b_E - D_{EI} w_I^*)$;

  \item $w_g = \tilde{C}_g^{-1}(b_g - D_{gE}w_E^*)$, $g=1,\cdots,G$;

  \item (the corrected electron segment) $w_E = w_E^* - C_E^{-1} \sum_{g=1}^G D_{Eg}w_g$;

  \item (the corrected ion segment) $w_I = w_I^* - A_I^{-1} D_{IE}w_E$.

\end{enumerate}
This preconditioning strategy is denoted as Schur1,
whose error can be formulated as follows:
\begin{align*}
  Err_{\textrm{Schur1}}(i,j) = \left \{
    \begin{aligned}
      & -D_{iE}C_E^{-1}D_{Ej},~~i\ne j,~i,j=1,\cdots,G \\
      & O, \qquad\qquad\quad~ \textrm{otherwise}
    \end{aligned}
  \right.,~~i,j=1,\cdots,G,E,I,
\end{align*}
where $O$ denotes the zero matrix of suitable size.


The last preconditioning strategy, denoted by Schur2,
is motivated by reducing the $(G+2)\times(G+2)$ block system into
a $2\times 2$ and $(G+1)\times(G+1)$ block systems.
It is of approximate block factorization type. Assembling the block decompositions
\begin{eqnarray*}
  \left[
  \begin{array}{cccc}
    A_1    &        &        & D_{1E} \\
           & \ddots &        & \vdots \\
           &        & A_G    & D_{GE} \\
    D_{E1} & \cdots & D_{EG} & A_E
  \end{array}\right]=\left[
  \begin{array}{ccc;{0.5pt/2pt}c}
        &        &        & D_{1E} \\
           & S_R &        & \vdots \\
           &        &     & D_{GE} \\ \hdashline[0.5pt/2pt]
     &  &  & A_E
  \end{array}\right]\left[
  \begin{array}{cccc}
    I    &        &        &  \\
           & \ddots &        &  \\
           &        & I    &  \\
    A_E^{-1}D_{E1} & \cdots & A_E^{-1}D_{EG} & I
  \end{array}\right]
\end{eqnarray*}
and
\begin{eqnarray*}
  \left[
  \begin{array}{cc}
    A_E    & D_{EI} \\
    D_{IE} & A_I
  \end{array}\right]=\left[
  \begin{array}{cc}
    A_E & O \\
    D_{IE} & S_I
  \end{array}\right]\left[
  \begin{array}{cc}
    I & A_E^{-1}D_{EI} \\
    O & I
  \end{array}\right],
\end{eqnarray*}
we obtain the full Schur complement preconditioning matrix of the form
\begin{eqnarray}\label{eq-3-02}
  \textrm{Schur2} & = & \left[
  \begin{array}{ccc;{0.5pt/2pt}cc}
        &        &       & D_{1E} & \\
           & S_R &       & \vdots & \\
           &        &    & D_{GE} & \\  \hdashline[0.5pt/2pt]
          &  &  & A_E & O \\
           &        & & D_{IE} & S_I
  \end{array}\right]\left[
  \begin{array}{ccccc}
    I    &        &       & &  \\
           & \ddots &     &   &  \\
           &        & I   & &  \\
    A_E^{-1}D_{E1} & \cdots & A_E^{-1}D_{EG} & I & A_E^{-1}D_{EI} \\
           &        &  & O & I
  \end{array}\right] \\ \nonumber & = & \left[
  \begin{array}{ccccc}
    A_1    &        &        & D_{1E} & D_{1E}A_E^{-1}D_{EI} \\
           & \ddots &        & \vdots & \vdots \\
           &        & A_G    & D_{GE} & D_{GE}A_E^{-1}D_{EI} \\
    D_{E1} & \cdots & D_{EG} & A_E & D_{EI} \\
  D_{IE}A_E^{-1}D_{E1} & \cdots & D_{IE}A_E^{-1}D_{EG} & D_{IE} & A_I
  \end{array}\right],
\end{eqnarray}
where $S_I$ and $S_R$ are the Schur complement matrices, respectively defined by
\begin{eqnarray*}
  S_I = A_I - D_{IE}A_E^{-1}D_{EI},\quad S_R = \left[ \begin{array}{ccc}
    A_1    &        &    \\
           & \ddots &    \\
           &        & A_G
  \end{array}\right]-\left[
  \begin{array}{c}
    D_{1E} \\
    \vdots \\
    D_{GE}
  \end{array}\right] A_E^{-1} \left[
  \begin{array}{ccc}
    D_{E1} & \cdots & D_{EG}
  \end{array}\right].
\end{eqnarray*}
The first term in \eqref{eq-3-02} represents the ``elimination" stage,
while the second factor in \eqref{eq-3-02} executes the ``substitution" stage.
Schur2 can be viewed as an extension to the matrix ${\bf A}$ of the primitive
variable block Schur complement preconditioner proposed recently by Weston et al. \cite{w-03}.
In the practical implementation, similar to Schur1, we take $\tilde{S}_R = {\bf diag}(\tilde{S}_1,\cdots,\tilde{S}_G)$
with $\tilde{S}_g=A_g-D_{gE}A_E^{-1}D_{Eg}$.
In this setting, the preconditioning operation proceeds as follows:
\begin{enumerate}

  \item (the intermediate electron segment) $w_E^* = A_E^{-1}b_E$;

  \item $w_g = \tilde{S}_g^{-1}(b_g-D_{gE}w_E^*)$, $g=1,\cdots,G$; \quad $w_I = S_I^{-1}(b_I-D_{IE}w_E^*)$;

  \item (the corrected electron segment) $w_E = w_E^* - A_E^{-1}(\sum_{g=1}^G D_{Eg}w_g + D_{EI}w_I)$.

\end{enumerate}
Observe that the action of Schur2 contains $\tilde{S}_g^{-1}$ ($g=1,\cdots,G$), $S_I^{-1}$ once and $A_E^{-1}$ twice.
The error associated with this preconditioner can be written as
\begin{align*}
  Err_{\textrm{Schur2}}(i,j) = \left \{
    \begin{aligned}
      & -D_{iE}A_E^{-1}D_{Ej},~~i\ne j,~i,j=1,\cdots,G,I \\
      & O, \qquad\qquad\quad~ \textrm{otherwise}
    \end{aligned}
  \right.,~~i,j=1,\cdots,G,E,I.
\end{align*}
By comparing with the block matrix $Err_{\textrm{Schur1}}$,
it is not difficult to see that there are $2G$ more nonzero data items in $Err_{\textrm{Schur2}}$.
A consequence of this is that, Schur1 may be able to perform better than Schur2 in terms of the convergence behavior.

\subsection{Fast approximations of the matrix inverse}\label{sec3.5}

It is customary to take advantage of the diagonal dominance property
in order to efficiently approximate the action of each matrix inverse of
the preceding block preconditioners. We start with a definition.

\begin{definition}\label{def2}
 For a given threshold $\theta_{wd}\in(0,1)$,
 define the weak diagonally dominant factor for the diagonal blocks $A_\alpha=(a_{kj}^{(\alpha)})_{N\times N}$ ($\alpha=1,\cdots,G,E,I$):
\begin{eqnarray*}
  \gamma_{wd}^{\alpha}(\theta_{wd})=\frac{\big{|}\{k:
  \sum_ja_{kj}^{(\alpha)}<\theta_{wd}\cdot a_{kk}^{(\alpha)},k=1,\cdots,N\}\big{|}}{N}\in[0,1].
\end{eqnarray*}
\end{definition}

\begin{remark}
 Besides $\gamma_{wd}^{\alpha}$ is a strictly increasing function,
 the larger the factor $\gamma_{wd}^{\alpha}(\theta_{wd})$,
 the weaker diagonally dominant the matrix $A_\alpha$ becomes.
\end{remark}

It is clear by Definition \ref{def2} that $A_\alpha$ can be regarded as a diagonal matrix
if the indicator $\gamma_{wd}^{\alpha}(\theta_{wd})=0$ holds.
$\theta_{wd}$ is taken to be $0.9$ in our experiments.
We implement two options for approximating the operation $w = A_\alpha^{-1}b$:
\begin{enumerate}

  \item[(\#1)] When $\gamma_{wd}^{\alpha}(\theta_{wd})=0$, this operation is performed by
  a fixed number of Jacobi iterations;

  \item[(\#2)] Otherwise, a fixed number of AMG V-cycles is used.

\end{enumerate}
It is worthwhile to point out that the interpolation operator
appeared in PCTL preconditioner should be constructed by another option:
\begin{enumerate}

  \item[(\#3)] The matrix inverse is solved by a number of AMG V-cycles to some prescribed accuracy.

\end{enumerate}
Also, notice that the inverse of $A_E$ still consists in the Schur complement matrices $\tilde{S}_R$ and $S_I$.
This can be avoided by a diagonal approximation. In this fashion, $\tilde{S}_R$ and $S_I$ are replaced by
\begin{eqnarray*}
  \hat{S}_R = {\bf diag}(\hat{S}_1,\cdots,\hat{S}_G)~~\mbox{with}~~
  \hat{S}_g = A_g - D_{gE}{\bf diag}(A_E)^{-1}D_{Eg},~g=1,\cdots,G;\quad
  \hat{S}_I = A_I - D_{IE}{\bf diag}(A_E)^{-1}D_{EI},
\end{eqnarray*}
where ${\bf diag}(A_E)$ denotes the diagonal matrix whose diagonal entries are those of $A_E$.
Therefore, $\hat{S}_g^{-1}$ ($g=1,\cdots,G$) and $\hat{S}_I^{-1}$
can be formed explicitly and approximated by option (\#3).
Approximations on $C_E$ and $\tilde{C}_g$ ($g=1,\cdots,G$),
denoted by $\hat{C}_E$ and $\hat{C}_g$,
are carried out in exactly the same manner. Alternatively,
another more accurate approach can be used to eliminate the quite slow convergence
when options (\#1)-(\#3) break down: the behavior of $w_I = S_I^{-1}q_I$
can be approximated in an iterative fashion, that is,
\begin{enumerate}

  \item[(\#4)] Taking the zero vector as the initial guess,
  the approximate solution at the $k$-th step is given by
  $$w_I^{(k+1)}=A_I^{-1}(q_I + D_{IE}A_E^{-1}D_{EI}w_I^{(k)}),$$
  where $A_g^{-1}$ ($g=E,I$) is iteratively approximated by option (\#1), (\#2) or (\#3).

\end{enumerate}
We see that option (\#4) would be much more expensive. In addition,
we can treat the remaining matrix inverses $C_E^{-1}$, $\tilde{C}_g$
and $\tilde{S}_g$ ($g=1,\cdots,G$) in the same fashion as above.

\subsection{An adaptive block preconditioning strategy}\label{sec3.6}

The coupling relationship can be used as the clear distinction to automatically choose a preconditioner
according to each different system from nonlinear iterations at each integration time-step.
Although the change from one linear system to the next may be relatively small,
the accumulative change after many time-steps would be significant.
An indicator is introduced in Definition \ref{def1} to measure the coupling strength for each coupling term.

\begin{definition}\label{def1}
 For a given threshold $\theta_{wc}\in(0,1)$,
 define the weak coupling factor for the coupling term
 $D_{\alpha\beta}=(d_{kj}^{(\alpha\beta)})_{N\times N}$ ($\alpha\ne\beta$)
 relative to $A_\alpha=(a_{kj}^{(\alpha)})_{N\times N}$:
\begin{eqnarray*}
  \gamma_{wc}^{\alpha\beta}(\theta_{wc})=\frac{\big{|}\{k:
  -d_{kk}^{(\alpha\beta)}\le\theta_{wc}\cdot a_{kk}^{(\alpha)},k=1,\cdots,N\}\big{|}}{N}\in[0,1].
\end{eqnarray*}
\end{definition}

\begin{remark}
 The larger the factor $\gamma_{wc}^{\alpha\beta}(\theta_{wc})$,
 the weaker coupling between the $\alpha$-th and $\beta$-th physical quantities appears.
 Furthermore, $\gamma_{wc}^{\alpha\beta}$ is a strictly increasing function.
\end{remark}

With this definition, it is possible to assume that $D_{\alpha\beta}$ can be neglected without a loss in robustness,
provided that $\gamma_{wc}^{\alpha\beta}(\theta_{wc})$ is larger than a switch criterion $\sigma_{wc}$.
The smaller $\sigma_{wc}$ is chosen, this is considered to be more aggressive.
In particular, if $\gamma_{wc}^{\alpha E}(\theta_{wc})>\sigma_{wc}$
is satisfied for some physical-variable $\alpha=1,\cdots,G,I$, then it can be extracted from the preconditioning matrix,
i.e., the $\alpha$-th row and column must be removed.
In this regard, the $\alpha$-th approximate solution is computed independently,
and the remaining unknowns can be treated in the same way as in the previous scenario.
In our computations, $\theta_{wc}$ and $\sigma_{wc}$ are set to $10^{-2}$ and $0.5$, respectively.
It should be mentioned that this adaptive block preconditioner strategy can be imposed on
the block preconditioners described in Sections \ref{sec3.3}-\ref{sec3.4}.

\begin{remark}
Obviously, the block lower (upper) triangular preconditioner is equivalent to the particular case when
$\gamma_{wc}^{EI}(\theta_{wc})>\sigma_{wc}$ and $\gamma_{wc}^{\alpha E}(\theta_{wc})>\sigma_{wc}$
($\gamma_{wc}^{IE}(\theta_{wc})>\sigma_{wc}$ and $\gamma_{wc}^{E\alpha}(\theta_{wc})>\sigma_{wc}$) for $\alpha=1,\cdots,G$,
while the block diagonal preconditioner corresponds to the simplest case that
$\gamma_{wc}^{\alpha E}(\theta_{wc})>\sigma_{wc}$ and $\gamma_{wc}^{E\alpha}(\theta_{wc})>\sigma_{wc}$ for $\alpha=1,\cdots,G,I$.
\end{remark}

\begin{remark}
 It must be emphasized that it is necessary to perform the coarse-grid correction step in PCTL preconditioner
 except when the conditions $\gamma_{wc}^{E\alpha}(\theta_{wc})>\sigma_{wc}$ ($\alpha=1,\cdots,G,I$) are all satisfied.
\end{remark}

\subsection{A summary of block preconditioning strategies}\label{sec3.7}

In summary, the block Schur complement preconditioners
utilize legacy algorithms to deal with the reduced scalar subsystems in an operator-splitting fashion,
whereas the PCTL preconditioner is derived from the inherence that there is no coupling
between the radiation and ion temperature variables, however,
they are separately coupled with the electron temperature variables.
Here, one should be aware that they are based on the \emph{divide-and-conquer} programming paradigm,
namely, each task is first split into more tractable subtasks,
and then their results are gathered and coupled.

Five types of operations that are performed in Setup and Preconditioning phases of
these block preconditioners are compared in Table \ref{ctp-03} with respect to a certain subsystem.
It is important to note that all the coupling terms $D_{\alpha\beta}$ ($\alpha\ne\beta$) are diagonal.
Consequently, all occurred matrix-vector multiplications reduce to
the so-called Hadamard products of two or three vectors except the one in the second step of PCTL preconditioner,
and the matrix update confined within diagonal entries is mathematically equivalent to one vector update.
Note that the additional $G+1$ matrix inverses as well as $G+2$ matrix updates are required in Setup phase of PCTL.
As can be seen, compared with Schur2,
there are $G+2$ more matrix-vector multiplications in PCTL,
as well as one more matrix inverse in Schur1.
Evidently, Schur2 is often preferred over PCTL and Schur1
when they exhibit similar inner and outer convergence behaviors.
It was mentioned above that Schur1 would probably be more accurate,
giving rise to a faster convergence and also better efficiency than Schur2 in many cases.

\begin{table}[htbp]
\tabcolsep=3pt
\footnotesize
\centering\caption{Operation comparisons in Setup and Preconditioning phases of PCTL, Schur1 and Schur2.}\label{ctp-03}\vskip 0.1cm
\begin{tabular}{ccccccccc}\hline
\multirow{2}{*}{}&\multicolumn{2}{c}{Setup phase}&&\multicolumn{4}{c}{Preconditioning phase} \\ \cline{2-3} \cline{5-9}
~ & matrix inverse & matrix update & & matrix inverse & matrix-vector multiplication & Hadamard product of vectors & vector update \\ \hline

 PCTL   & $G+1$ & $G+2$ &  & $G+3$ & $G+2$ & $5G+7$ & $6G+8$ \\
 Schur1 &   0   &   0   &  & $G+4$ &   0   & $3G+3$ & $2G+4$ \\ 
 Schur2 &   0   &   0   &  & $G+3$ &   0   & $3G+3$ & $2G+4$ \\ \hline

\end{tabular}
\end{table}

Each adaptive block preconditioned GMRES($m$) solver consists of two levels.
The adaptive block preconditioning strategy described in Section \ref{sec3.6} is considered at the first level.
Its purpose is to find an appropriate preconditioner with less computation time for the problems relatively easier to solve.
The second level bifurcates every matrix inverse appeared in
the preceding block preconditioners into the fast or accurate approximation proposed in Section \ref{sec3.5}.

\subsection{Implementation details}\label{sec3.8}

\begin{figure}[htbp]
\centerline{
\includegraphics[scale=0.3]{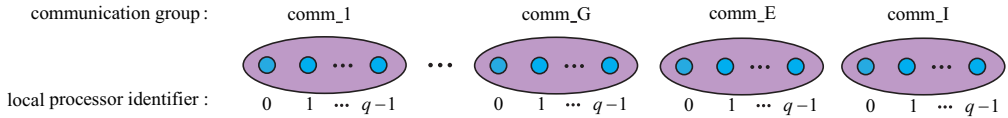}}
\caption{The communicator splitting and processor partitioning,
where $(G+2)q$ is the number of processors.}\label{fig-02}
\end{figure}

The three block preconditioners described in Section \ref{sec3.3}-\ref{sec3.4}
are implemented on the basis of the HYPRE package
using C and message passing interface (MPI) with data located in distributed memory.
The parallel executable code decomposes the global communicator, 
illustrated in Figure \ref{fig-02},
such that each subsystem owns a group holding an ordered collection of processor identifiers
for the purpose of running parallelized modules to solve large problems.

An effective parallel data decomposition plays a key role in achieving good performance.
Each subsystem should be partitioned in such a way that each part is of about equal size to avoid load balance issues.
In our parallel implementation, all diagonal blocks $A_\alpha$ are stored using the ParCSR matrix data structure
within the associated communication group,
while all coupling terms $D_{\alpha\beta}$ are represented by the ParVector data structure of global length $N$.
More detailed information about ParCSR and ParVector can be found in Falgout et al. \cite{f-01}.
Each matrix inverse is accomplished by BoomerAMG \cite{h-01} with one of options depicted in Section \ref{sec3.5},
where we use V(1,1)-cycles with hybrid Gauss-Seidel relaxation in symmetric ordering,
coarse-grid matrices formed algebraically by Galerkin process,
and algebraic interface-based coarsening recently presented by Xu and Mo in \cite{x-02}
to gain a good balance between convergence and complexity.

\begin{figure}[htbp]
\centerline{
\includegraphics[scale=0.25]{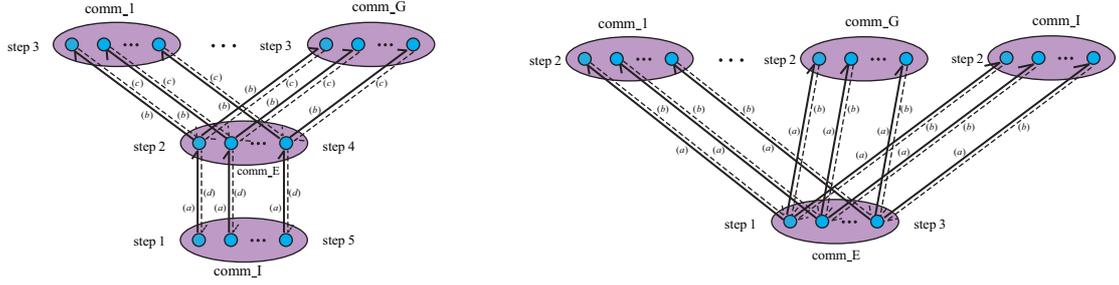}}
\caption{The parallel flow diagram and MPI communication pattern of Schur1 (left) and Schur2 (right).}\label{fig-03}
\end{figure}

It is apparent that Schur1 and Schur2 are far easier to implement than PCTL
because of the construction and numerical solution of the coarse-grid system,
involving much more complicated global communication between different processors \cite{y-04},
e.g., all diagonal blocks $A_\alpha$ ($\alpha=1,\cdots,G,E,I$)
and coupling terms $D_{\alpha\beta}$ ($\alpha\ne\beta$)
should be distributed locally (data sets used in setup and preconditioning phases)
and globally (data sets used to generate the globally distributed matrix $A_c$ by formula \eqref{eq-3-03} in setup phase).
However, their distributed local data sets are enough for Schur1 and Schur2,
whose parallel flow diagrams and MPI communication patterns are depicted in Figure \ref{fig-03}.
The numerical solution of a certain subsystem at step 1, 2, 3 and 4
is followed by the data exchange process labeled (a), (b), (c) and (d), respectively.
Another communication group is to be established in setup phase so that processors
with the same local identifier in the partition (shown in Figure \ref{fig-02}) are grouped together
for the purpose of making exchanging messages (such as (b) and (c) in Schur1) easier and more efficient to implement.
At this point, the two MPI functions {\em MPI\_Bcast} and {\em MPI\_Gather}
can be utilized to spread elements from the root processor (belonging to {\em comm\_E}) to the others and
take elements from other processors and gather them to the root processor, respectively.

\section{Numerical results}\label{sec4}

In this section, we provide experimental results intended to study
the effectiveness and parallel scalability of the block preconditioners against 
BoomerAMG \cite{h-01} implemented in the parallel C library HYPRE version 2.15.1.
These results are obtained on the Tianhe-1 supercomputer \cite{y-05},
which is a multi-array, configurable and cooperative parallel system with a theoretical peak speed of 1.372 petaflops,
composed of high performance general-purpose microprocessors and a high-speed Infiniband network.
Our code is compiled with mpich-3.1.3 using the icc compiler version 11.1.059 and -O2 optimization level.
The convergence of the preconditioned GMRES($m$) was halted when the Euclidean norm of
the current relative residual was smaller than $10^{-7}$.
Unless otherwise stated, the tolerances for stopping are chosen to be $10^{-2}$ for performance consideration
in approximating $\hat{C}_E^{-1}$, $\hat{S}_I^{-1}$, $\hat{C}_g^{-1}$ and $\hat{S}_g^{-1}$ ($g=1,\cdots,G$)
in Schur complement preconditioners,
as well as the construction of the interpolation and block relaxations in the PCTL preconditioner.

The following notations are given in tables below to illustrate results:
$N_{it}$ denotes the number of preconditioned GMRES($m$) iterations;
$T_{cpu}$ represents the CPU time measured in seconds to solve the linear system by the preconditioned GMRES($m$) solver;
$\alpha$PCTL, $\alpha$Schur1 and $\alpha$Schur2 stand for the adaptive PCTL, Schur1 and Schur2 variants, respectively;
$np$ indicates the number of used processors.

\subsection{Numerical experiments on one processor}

Results in this subsection have been structured in two distinct parts.
They are based on two suites of representative MGD linear systems from two-dimensional capsule implosion simulations:
the first suite includes 4 three-temperature (i.e., only one-group) linear systems,
respectively denoted by $S_1$-$S_4$; 
the other suite consists of 37 three-temperature and 20 twenty-group linear systems.

The performance of seven preconditioned GMRES(30) solvers is investigated for problems $S_1$-$S_4$
on two different grids, as shown in Table \ref{ctp-01}. We can observe from this table that
(i) AMG preconditioned GMRES(30) doesn't converge robustly enough,
and achieves convergence in an excessive number of iterations (shown by $S_2$);
(ii) PCTL, Schur1 and Schur2 exhibit a better convergence behavior, and,
as expected, Schur1 is numerically more robust than the other two preconditioners;
(iii) Schur1 and Schur2 require less wall time than that of PCTL;
(iv) The adaptive variants are much more efficient,
although there are additional arithmetic operations to determine an appropriate and inexpensive preconditioner.

\begin{table}[htbp]
\tabcolsep=1.8pt
\footnotesize
\centering\caption{Number of iterations and wall time of seven preconditioned GMRES(30) solvers.}\label{ctp-01}\vskip 0.1cm
\begin{tabular}{cccccccccccccccccccccccccccccc}\hline
\multirow{3}{*}{}&\multicolumn{14}{c}{$32,000\times24$}&\multirow{2}{*}{}&\multicolumn{14}{c}{$64,000\times48$} \\ \cline{2-15} \cline{17-30}
~&\multicolumn{2}{c}{AMG}&\multicolumn{2}{c}{PCTL}&\multicolumn{2}{c}{$\alpha$PCTL}
&\multicolumn{2}{c}{Schur1}&\multicolumn{2}{c}{$\alpha$Schur1}&\multicolumn{2}{c}{Schur2}&\multicolumn{2}{c}{$\alpha$Schur2}&\multirow{2}{*}{}
&\multicolumn{2}{c}{AMG}&\multicolumn{2}{c}{PCTL}&\multicolumn{2}{c}{$\alpha$PCTL}
&\multicolumn{2}{c}{Schur1}&\multicolumn{2}{c}{$\alpha$Schur1}&\multicolumn{2}{c}{Schur2}&\multicolumn{2}{c}{$\alpha$Schur2} \\ \cline{2-15} \cline{17-30}
~& $N_{it}$ & $T_{cpu}$ & $N_{it}$ & $T_{cpu}$ & $N_{it}$ & $T_{cpu}$ & $N_{it}$ & $T_{cpu}$ & $N_{it}$ & $T_{cpu}$ & $N_{it}$ & $T_{cpu}$ & $N_{it}$ & $T_{cpu}$ &
& $N_{it}$ & $T_{cpu}$ & $N_{it}$ & $T_{cpu}$ & $N_{it}$ & $T_{cpu}$ & $N_{it}$ & $T_{cpu}$ & $N_{it}$ & $T_{cpu}$ & $N_{it}$ & $T_{cpu}$ & $N_{it}$ & $T_{cpu}$ \\ \hline

 $S_1$ &  9 &  3.6 &  3 &  1.7 &  2 & 1.5 & 3 & 1.4 & 3 & 1.4 & 4 & 1.6 & 4 & 1.6 &  &   9 &  14.8 &  3 &   7.6 &  2 &  6.4 & 3 &  7.1 & 3 &  6.9 & 5 &  7.7 & 5 &  7.8 \\ 
 $S_2$ & 80 & 31.2 & 12 &  6.1 & 10 & 5.3 & 4 & 3.1 & 5 & 2.9 & 4 & 3.0 & 5 & 2.9 &  & 100 & 164.2 & 17 & 270.3 & 14 & 28.7 & 4 & 13.5 & 6 & 14.0 & 4 & 11.8 & 6 & 13.0 \\ 
 $S_3$ & 29 & 12.1 & 11 &  5.4 &  6 & 3.8 & 4 & 3.1 & 5 & 2.9 & 6 & 3.4 & 6 & 3.2 &  &  32 &  56.4 & 14 &  28.2 &  7 & 17.4 & 4 & 14.3 & 6 & 13.6 & 8 & 18.0 & 8 & 16.2 \\ 
 $S_4$ & 11 &  4.8 &  3 & 11.6 &  6 & 3.9 & 3 & 9.1 & 6 & 3.9 & 3 & 9.3 & 6 & 3.9 &  &  11 &  20.6 &  3 &  24.9 &  6 & 16.3 & 3 & 18.5 & 6 & 16.3 & 3 & 19.5 & 6 & 16.3 \\ \hline

\end{tabular}
\end{table}

An important observation to make in Table \ref{ctp-01} is why
the computational cost of PCTL becomes unacceptable for $S_2$ on $64,000\times48$ grid,
despite its reasonable iteration counts.
The focus is on the inner iteration comparisons among PCTL, $\alpha$PCTL, Schur1, $\alpha$Schur1, Schur2 and $\alpha$Schur2.
The results are illustrated in Table \ref{ctp-02},
where (i) entries with the superscript {\footnotesize $*$} indicate that option (\#1) is used and option (\#2) otherwise,
here we consider just a single Jacobi iteration or AMG V-cycle to approximate matrix inverse;
(ii) entries with the superscript {\footnotesize $\dagger$} inform that
the block diagonal preconditioner with option (\#1) or (\#2) is utilized to improve performance.
Since the simulation is much more sensitive to the degree of accuracy in radiation temperatures,
three iterations are used to solve approximately for them, while one iteration is applied to subsystems
associated with electron and ion temperatures.
The primary reason why PCTL is rather expensive is that a total number of 623 AMG V-cycles are required to solve all of the subsystems
associated with all radiation temperatures. This is remedied in $\alpha$PCTL by using option (\#1) or (\#2)
instead of option (\#3), however, with many more coarse-grid correction steps (from 9 to 15).
Fortunately, there is a decrease in the number of outer iterations (from 17 to 14),
resulting in a much faster preconditioner.

\begin{table}[htbp]
\tabcolsep=2pt
\footnotesize
\centering\caption{Inner iteration comparisons of six block preconditioners for the $64,000\times48$ grid.}\label{ctp-02}\vskip 0.1cm
\begin{tabular}{ccccccccccccccccccccccccccccccccc}\hline
\multirow{2}{*}{}&\multicolumn{6}{c}{PCTL}&&\multicolumn{6}{c}{$\alpha$PCTL}
&&\multicolumn{3}{c}{Schur1}&&\multicolumn{3}{c}{$\alpha$Schur1}
&&\multicolumn{3}{c}{Schur2}&&\multicolumn{3}{c}{$\alpha$Schur2} \\ \cline{2-7} \cline{9-14} \cline{16-18} \cline{20-22} \cline{24-26} \cline{28-30}
~& $p_1$ & $p_I$ & $A_E^{-1}$ & $A_1^{-1}$ & $A_I^{-1}$ & $A_c^{-1}$ && $p_1$ & $p_I$ & $A_E^{-1}$ & $A_1^{-1}$ & $A_I^{-1}$ & $A_c^{-1}$
&& $C_E^{-1}$ & $\hat{C}_1^{-1}$ & $A_I^{-1}$ && $C_E^{-1}$ & $\hat{C}_1^{-1}$ & $A_I^{-1}$
&& $A_E^{-1}$ & $\hat{S}_1^{-1}$ & $\hat{S}_I^{-1}$ && $A_E^{-1}$ & $\hat{S}_1^{-1}$ & $\hat{S}_I^{-1}$ \\ \hline

 $S_1$ & 0 & 1 &  4 &   3 &  3 &  4 &  & 1 & 1 & $3^*$ &  3 &  $3^*$ &  3 &  & 6 & 3 & 7 &  &  6 & 3 &  $7^*$ &  & 14 &  3 & 5 &  & $11^*$ & 4 & $5^*$ \\ 
 $S_2$ & 2 & 1 & 10 & 623 & 9 & 9 &  & 6 & 1 &    15 & 15 & $15^*$ & 15 &  & 10 & 12 & 9 &  & 14 & 7 & $13^*$ &  &  9 & 10 & 5 &  & 13 & 7 & $7^*$ \\ 
 $S_3$ & 3 & 1 &  5 &  19 & 5 & 11 &  & 6 & 1 & 8 & 8 & $8^*$ & 8 &  & 10 & 14 & 10 &  & 14 & 7 & $12^*$ &  & 17 & 14 & 9 && 17 & 9 & $9^*$ \\ 
 $S_4$ & 7 & 1 &  6 &  27 & 4 & 1 &  & 0 & 0 & $7^\dagger$ & $21^\dagger$ & $7^\dagger$ &  0
 &  & 5 & 25 & 5 &  & $7^\dagger$ & $21^\dagger$ & $7^\dagger$ & &  7 & 27 & 4 && $7^\dagger$ & $21^\dagger$ & $7^\dagger$ \\ \hline

\end{tabular}
\end{table}

\begin{figure}[htbp]
\centerline{
\includegraphics[scale=0.45]{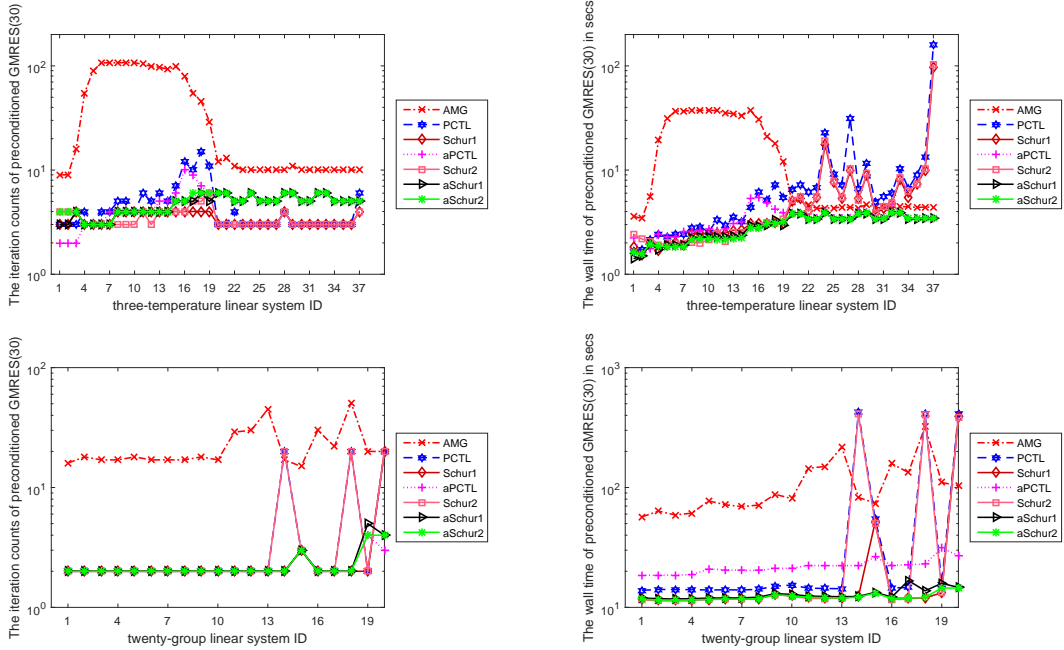}}
\caption{The iteration counts and wall time 
as functions of representative three-temperature (top) and twenty-group (bottom) linear systems.}\label{fig-01}
\end{figure}

Next we examine the iteration counts and wall time as functions of the second suite of representative linear systems.
A plot of this data appears in Figure \ref{fig-01} with problem sizes roughly 2.3 and 16.9 million DoFs, respectively.
It can be easily seen that PCTL, Schur1, Schur2 and their adaptive variants converge much faster than AMG, however,
with 3 twenty-group exceptions which could not be solved by PCTL, Schur1 and Schur2.
The reason for this is that the tolerance for stopping in option (\#3) is not small enough.
However, the smaller the tolerance, the higher their total costs.
It should be emphasized that, although they may demonstrate a slightly slower convergence,
$\alpha$PCTL, $\alpha$Schur1 and $\alpha$Schur2 are less expensive to use,
because the actual convergence behavior is produced by the block diagonal preconditioner.
From the viewpoint of the total wall time,
$\alpha$Schur1 and $\alpha$Schur2 have comparable computational cost, 
and they are both observed to be the most computationally efficient option,
resulting in a 13.7\% (three-temperature) and 41.6\% (twenty-group) reduction compared with $\alpha$PCTL,
as well as 5.5 (three-temperature) and 8.5 (twenty-group) times faster than AMG.



\subsection{Parallel experiments}

Run times, strong and weak parallel scaling properties are both of particular interest in this subsection.
Since our previous performance evaluations showed that $\alpha$PCTL, $\alpha$Schur1 and $\alpha$Schur2
leads generally to better convergence, robustness and efficiency,
we only coded their parallel implementations while PCTL, Schur1 and Schur2 are not considered here.
Timings are measured by the MPI function {\em MPI\_Wtime} in seconds.

We consider 6 twenty-group linear systems, respectively denoted by $M_1$-$M_6$,
from the previous subsection and examine the strong parallel scalabilities of AMG,
$\alpha$PCTL, $\alpha$Schur1 and $\alpha$Schur2 preconditioned GMRES(30) solvers
on a $16,000\times48$ grid (see Table \ref{ctp-04}).
The simulation is run on 22 to 352 cores, each time doubling the number of used cores.
The problem size on each processor is 48,000 DoFs for 352 cores.
The results in Table \ref{ctp-04} show that these four preconditioners exhibit
good convergence (in terms of number of iterations) and numerical (regarding wall time) scalabilities in the strong sense.
We note that all of $\alpha$PCTL, $\alpha$Schur1 and $\alpha$Schur2 preconditioned GMRES(30) solvers converge in at most 3 steps
robustly with respect to the number of processors and problem character,
and run averagely 3.4, 4.1 and 5.5 times faster than AMG on 352 cores.

\begin{table}[htbp]
\tabcolsep=2pt
\footnotesize
\centering\caption{Number of iterations and wall time
of four preconditioned GMRES(30) solvers on a $16,000\times48$ grid
in a strong scaling study.}\label{ctp-04}\vskip 0.1cm
\begin{tabular}{ccccccccccccccccccccccccccccccccccccccc}\hline
\multirow{3}{*}{}&\multicolumn{10}{c}{AMG}&&\multicolumn{10}{c}{$\alpha$PCTL} \\ \cline{2-11} \cline{13-22}
$np$& \multicolumn{2}{c}{$22$} & \multicolumn{2}{c}{$44$} & \multicolumn{2}{c}{$88$} & \multicolumn{2}{c}{$176$} & \multicolumn{2}{c}{$352$}
&& \multicolumn{2}{c}{$22$} & \multicolumn{2}{c}{$44$} & \multicolumn{2}{c}{$88$} & \multicolumn{2}{c}{$176$} & \multicolumn{2}{c}{$352$} \\ \cline{2-11} \cline{13-22}
~&$N_{it}$&$T_{cpu}$&$N_{it}$&$T_{cpu}$&$N_{it}$&$T_{cpu}$&$N_{it}$&$T_{cpu}$&$N_{it}$&$T_{cpu}$
&&$N_{it}$&$T_{cpu}$&$N_{it}$&$T_{cpu}$&$N_{it}$&$T_{cpu}$&$N_{it}$&$T_{cpu}$&$N_{it}$&$T_{cpu}$ \\ \hline

 $M_1$ & 19 & 22.7 & 18 & 14.7 & 17 &  9.4 & 19 & 5.8 & 16 & 3.0 && 2 & 6.8 & 2 & 4.3 & 2 & 2.5 & 2 & 1.7 & 2 & 1.1 \\ 
 $M_2$ & 18 & 22.0 & 18 & 15.0 & 18 &  9.9 & 18 & 5.6 & 17 & 3.1 && 2 & 6.7 & 2 & 4.3 & 2 & 2.5 & 2 & 1.8 & 2 & 1.2 \\ 
 $M_3$ & 18 & 21.9 & 17 & 14.5 & 19 & 10.2 & 19 & 5.8 & 19 & 3.3 && 2 & 6.8 & 2 & 4.2 & 2 & 2.5 & 2 & 1.7 & 2 & 1.1 \\ 
 $M_4$ & 34 & 55.6 & 32 & 23.1 & 36 & 16.5 & 32 & 8.2 & 34 & 4.5 && 2 & 6.9 & 2 & 4.2 & 2 & 2.6 & 2 & 1.6 & 2 & 1.1 \\ 
 $M_5$ & 26 & 53.7 & 24 & 22.5 & 24 & 13.7 & 27 & 7.4 & 30 & 4.3 && 2 & 7.0 & 2 & 4.3 & 2 & 2.6 & 2 & 1.6 & 2 & 1.1 \\ 
 $M_6$ & 27 & 49.8 & 29 & 31.3 & 33 & 17.8 & 29 & 7.6 & 32 & 4.4 && 3 & 8.5 & 3 & 5.2 & 2 & 2.7 & 3 & 2.1 & 2 & 1.1 \\ 
\\
\multirow{3}{*}{}&&\multicolumn{10}{c}{$\alpha$Schur1}&&\multicolumn{10}{c}{$\alpha$Schur2} \\ \cline{2-11} \cline{13-22}
$np$& \multicolumn{2}{c}{$22$} & \multicolumn{2}{c}{$44$} & \multicolumn{2}{c}{$88$} & \multicolumn{2}{c}{$176$} & \multicolumn{2}{c}{$352$}
&& \multicolumn{2}{c}{$22$} & \multicolumn{2}{c}{$44$} & \multicolumn{2}{c}{$88$} & \multicolumn{2}{c}{$176$} & \multicolumn{2}{c}{$352$} \\ \cline{2-11} \cline{13-22}
~&$N_{it}$&$T_{cpu}$&$N_{it}$&$T_{cpu}$&$N_{it}$&$T_{cpu}$&$N_{it}$&$T_{cpu}$&$N_{it}$&$T_{cpu}$
&&$N_{it}$&$T_{cpu}$&$N_{it}$&$T_{cpu}$&$N_{it}$&$T_{cpu}$&$N_{it}$&$T_{cpu}$&$N_{it}$&$T_{cpu}$ \\ \hline

 $M_1$ & 2 & 6.8 & 2 & 4.0 & 2 & 2.5 & 2 & 1.6 & 2 & 0.9 && 2 & 6.1 & 2 & 3.8 & 2 & 2.1 & 2 & 1.3 & 2 & 0.7 \\ 
 $M_2$ & 2 & 6.7 & 2 & 4.1 & 2 & 2.5 & 2 & 1.5 & 2 & 0.9 && 2 & 6.1 & 2 & 3.9 & 2 & 2.1 & 2 & 1.2 & 2 & 0.7 \\ 
 $M_3$ & 2 & 6.8 & 2 & 4.1 & 2 & 2.4 & 2 & 1.4 & 2 & 0.9 && 2 & 6.2 & 2 & 3.8 & 2 & 2.3 & 2 & 1.1 & 2 & 0.6 \\ 
 $M_4$ & 2 & 6.8 & 2 & 4.2 & 2 & 2.5 & 2 & 1.5 & 2 & 0.9 && 2 & 6.2 & 2 & 3.9 & 2 & 2.2 & 2 & 1.3 & 2 & 0.7 \\ 
 $M_5$ & 2 & 6.8 & 2 & 4.1 & 2 & 2.5 & 3 & 1.9 & 2 & 1.0 && 2 & 6.2 & 2 & 4.0 & 2 & 2.3 & 2 & 1.2 & 2 & 0.7 \\ 
 $M_6$ & 3 & 8.2 & 3 & 5.0 & 2 & 2.5 & 3 & 1.9 & 2 & 0.9 && 3 & 7.4 & 3 & 4.7 & 2 & 2.2 & 3 & 1.6 & 2 & 0.7 \\ \hline

\end{tabular}
\end{table}

The subsequent experiments are run for a weak parallel scaling study
to analyze the convergence and numerical scalabilities 
of GMRES(30) solvers preconditioned by AMG, $\alpha$PCTL, $\alpha$Schur1 and $\alpha$Schur2.
We start with a mesh size of $4,000\times12$ and then refine the mesh in all directions up to $16,000\times48$,
so the largest problem has about 16.9 million unknowns.
We run these problems using 22, 88 and 352 cores (with 48,000 DoFs per processor), respectively.
As shown in Table \ref{ctp-05}, all of the preconditioners tested here weakly scale well in the numerical sense.
$\alpha$PCTL, $\alpha$Schur1 and $\alpha$Schur2 preconditioned GMRES(30) solvers converge robustly with respect to the discrete problem size,
while the number of iterations of AMG preconditioned GMRES(30) solver for $M_5$ varies from
16 on the $4,000\times12$ grid to 30 on the $16,000\times48$ grid to achieve convergence.
Regarding the parallel efficiency calculated by $T_{cpu}^{(88)} / T_{cpu}^{(352)}$ \cite{s-05},
where $T_{cpu}^{(k)}$ is the wall time using $k$ processors,
the average parallel efficiency of AMG, $\alpha$PCTL, $\alpha$Schur1 and $\alpha$Schur2 preconditioned GMRES(30) solvers
is 72.2\%, 79.3\%, 74.7\% and 75.1\%, respectively. Furthermore, in terms of wall time,
$\alpha$Schur2 has obvious advantages over AMG, $\alpha$PCTL and $\alpha$Schur1.

\begin{table}[htbp]
\tabcolsep=2pt
\footnotesize
\centering\caption{Number of iterations and wall time
of four preconditioned GMRES(30) solvers
in a weak scaling study.}\label{ctp-05}\vskip 0.1cm
\begin{tabular}{ccccccccccccccccccccccccccccccccccccccc}\hline
\multirow{3}{*}{}&\multicolumn{6}{c}{AMG}&&\multicolumn{6}{c}{$\alpha$PCTL}
&&\multicolumn{6}{c}{$\alpha$Schur1}&&\multicolumn{6}{c}{$\alpha$Schur2} \\ \cline{2-7} \cline{9-14} \cline{16-21} \cline{23-28}
$np$& \multicolumn{2}{c}{$22$} & \multicolumn{2}{c}{$88$} & \multicolumn{2}{c}{$352$}
&& \multicolumn{2}{c}{$22$} & \multicolumn{2}{c}{$88$} & \multicolumn{2}{c}{$352$}
&& \multicolumn{2}{c}{$22$} & \multicolumn{2}{c}{$88$} & \multicolumn{2}{c}{$352$}
&& \multicolumn{2}{c}{$22$} & \multicolumn{2}{c}{$88$} & \multicolumn{2}{c}{$352$}
\\ \cline{2-7} \cline{9-14} \cline{16-21} \cline{23-28}
~&$N_{it}$&$T_{cpu}$&$N_{it}$&$T_{cpu}$&$N_{it}$&$T_{cpu}$
&&$N_{it}$&$T_{cpu}$&$N_{it}$&$T_{cpu}$&$N_{it}$&$T_{cpu}$
&&$N_{it}$&$T_{cpu}$&$N_{it}$&$T_{cpu}$&$N_{it}$&$T_{cpu}$
&&$N_{it}$&$T_{cpu}$&$N_{it}$&$T_{cpu}$&$N_{it}$&$T_{cpu}$ \\ \hline

 $M_1$ & 15 & 1.33 & 15 & 2.29 & 16 & 3.01 && 2 & 0.62 & 2 & 0.87 & 2 & 1.08 && 2 & 0.48 & 2 & 0.67 & 2 & 0.91 && 2 & 0.35 & 2 & 0.50 & 2 & 0.71 \\ 
 $M_2$ & 14 & 1.25 & 17 & 2.49 & 17 & 3.13 && 2 & 0.58 & 2 & 0.85 & 2 & 1.15 && 2 & 0.46 & 2 & 0.65 & 2 & 0.93 && 2 & 0.36 & 2 & 0.49 & 2 & 0.69 \\ 
 $M_3$ & 14 & 1.24 & 18 & 2.63 & 19 & 3.26 && 2 & 0.58 & 2 & 0.84 & 2 & 1.11 && 2 & 0.47 & 2 & 0.65 & 2 & 0.91 && 2 & 0.35 & 2 & 0.48 & 2 & 0.64 \\ 
 $M_4$ & 28 & 2.05 & 30 & 3.55 & 34 & 4.52 && 2 & 0.59 & 2 & 0.88 & 2 & 1.12 && 2 & 0.46 & 2 & 0.64 & 2 & 0.89 && 2 & 0.35 & 2 & 0.46 & 2 & 0.67 \\ 
 $M_5$ & 16 & 1.47 & 18 & 2.61 & 30 & 4.31 && 2 & 0.58 & 3 & 0.98 & 2 & 1.11 && 2 & 0.45 & 3 & 0.76 & 2 & 0.95 && 2 & 0.34 & 3 & 0.58 & 2 & 0.70 \\ 
 $M_6$ & 21 & 1.68 & 21 & 2.78 & 32 & 4.40 && 2 & 0.59 & 2 & 0.86 & 2 & 1.09 && 2 & 0.46 & 3 & 0.76 & 2 & 0.94 && 2 & 0.34 & 3 & 0.59 & 2 & 0.72 \\ \hline

\end{tabular}
\end{table}

\section{Conclusions}\label{sec5}

Three types of AMG block-based preconditioning techniques and two improvements have been introduced
to solve the linear systems resulting from 
MGD equations. They are easy to set up, due to the fact that they only involve operations on blocks
that are readily extracted from the monolithic linear system.
Our experiments show that they all scale well both algorithmically and in parallel,
and exhibit overall better convergence behavior and computational efficiency than the monolithic AMG preconditioner.

\section*{Acknowledgement}

This work is financially supported by Science Challenge Project (TZ2016002),
National Natural Science Foundation of China (11601462, 11971414),
Project of Scientific Research Fund of Hunan Provincial Science and Technology Department (2018WK4006),
Hunan Provincial Civil-Military Integration Industrial Development Project ``Adaptive Multilevel Solver and Its Application in ICF Numerical Simulation"
and Hunan Provincial Natural Science Foundation of China (2018JJ3494).
The numerical calculations in this paper have been done on the supercomputing system of the National Supercomputing Center in Changsha.

\end{document}